\documentclass{amsart}
\usepackage{amscd,amssymb}

\theoremstyle{plain}
\newtheorem{thm}{Theorem}[section]
\newtheorem{lem}[thm]{Lemma}

\newtheorem{prop}[thm]{Proposition}

\theoremstyle{definition}

\theoremstyle{remark}

\numberwithin{equation}{section}

\begin{document}
\setcounter{section}{-1}

Aleksander Ivanov
\footnote{University of Wroc{\l}aw, Institute of Mathematics,
pl.Grunwaldzki 2/4, 50-384 Wroc{\l}aw, Poland. \\
e-mail: ivanov@math.uni.wroc.pl}
\footnote{The research was supported by KBN grants 2 P03A 007 19
and 1 P03A 025 28}
\bigskip

\centerline{{\bf AN $\aleph_{0}$-CATEGORICAL THEORY WHICH IS
NOT G-COMPACT}}
\centerline{{\bf AND DOES NOT HAVE AZ-ENUMERATIONS}}

\bigskip


{\bf Abstract.}
\footnote{{\em Key words or phrases:} $\aleph_{0}$-categorical theories,
G-compactness, AZ-enumerations}
We give an example of an $\aleph_{0}$-categorical theory
which is not G-compact.
The countable model of this theory does not have AZ-enumerations.

\bigskip

\section{Introduction}
G-Compact theories were introduced in \cite{lascar} in 1982.
The first examples of non-G-compact theories were found in 2000
(see \cite{clpz}).
The main result of the paper is a construction
of an $\aleph_{0}$-categorical theory which is not G-compact.
We have also found that the countable model of this theory
does not have AZ-enumerations.
This answers a question formulated in \cite{AlCh}.
\bigskip

Let $T$ be a first-order theory over a countable language.
It is assumed that models of $T$ are elementary substructures
of a sufficiently saturated monster model $\mathbb{C}$.
We use $A,B,C$ to denote subsets of $\mathbb{C}$, assumed
to be much smaller than $\mathbb{C}$.
If $\bar{a}$ is a tuple from a model $M$, we often abuse
notation by writing $\bar{a}\in M$.
If $r(\bar{x})$ is a type in (possibly infinitely many )
variables $\bar{x}$, we denote by $r(M)$ the set of tuples
(sequences) from $M$ which realize $r$.
For any structure $M$ and $A\subseteq M$, define $Aut(M/A)$ to
be the group of automorphisms of $M$ which fix $A$ pointwise.
\parskip0pt

The following definitions and facts are partially taken from
\cite{clpz} and \cite{laspil}.
For $\delta\in\{ 1,2,...,\omega\}$ and a set $A$
let $E_{L}^{A,\delta}$ be the finest bounded
$Aut(\mathbb{C}/A)$-invariant equivalence relation on
$\delta$-sequences.
The classes of $E_{L}^{A,\delta}$ are called Lascar strong
types.
The relation $E^{A,\delta}_{L}$ can be characterized as follows:
$(\bar{a},\bar{b})\in E^{A,\delta}_{L}$ if there are models
$M_{1},...,M_{n}<\mathbb{C}$ containing $A$, and sequences
$\bar{a}_{0}(=\bar{a}),...,\bar{a}_{n}(=\bar{b})$ such that
$tp(\bar{a}_{i}/M_{i+1})=tp(\bar{a}_{i+1}/M_{i+1})$, $0\le i<n$.
Equivalently $(\bar{a},\bar{b})\in E^{A,\delta}_{L}$
if there are sequences
$\bar{a}_{0}(=\bar{a}),...,\bar{a}_{n}(=\bar{b})$ such that
each pair $\bar{a}_{i},\bar{a}_{i+1}$, $0\le i<n$, extends to
an infinite indiscernible sequence over $A$.\parskip0pt

Let $E^{A,\delta}_{KP}$ be the finest bounded $A$-type-definable
equivalence relation on $\delta$-sequences and $E^{A,\delta}_{Sh}$
be the intersection of all finite $A$-definable equivalence
relations on $\delta$-sequences.
Sequences $\bar{a}$ and $\bar{b}\in X$ have the same strong
(KP-strong) types over $A$ if and only if they are
$E^{A,\delta}_{Sh}$-equivalent ($E^{A,\delta}_{KP}$-equivalent).
It is known from \cite{lascar} that for $\omega$-categorical
theories and for finite $A$ and $\delta$,
$E^{A,\delta}_{Sh}=$ $E^{A,\delta}_{KP}=$ $E^{A,\delta}_{L}$.
\parskip0pt

Let $M$ be a saturated structure of uncountable cardinality
and let $Aut_{L}(M)$ be the group of all Lascar strong
automorphisms (fixing the classes of all bounded invariant
equivalence relations).
Then $Gal_{L}(Th(M))=Aut(M)/Aut_{L}(M)$,
the {\em Galois group of} $Th(M)$, does not depend on $M$.
Let $Aut_{KP}(M)$ be the subgroup of $Aut(M)$ consisting
of all automorphisms fixing the classes of all bounded
$\emptyset$-type-definable equivalence relations.
It makes sense to consider the following group extension:
$$
1 \rightarrow Aut_{KP}(M)/Aut_{L}(M) \rightarrow
Gal_{L} \rightarrow Gal_{KP}=Aut(M)/Aut_{KP}(M) \rightarrow 1.
$$
It is known that the group $Gal_{KP}$ is compact.
The theory $Th(M)$ is called {\em G-compact} if
$Gal_{L}=Gal_{KP}$.
It is equivalent to either of the following conditions:
(1) $E^{\emptyset ,\delta}_{L}= E^{\emptyset ,\delta}_{KP}$ 
for all $\delta$,
or (2) $E^{\emptyset ,\delta}_{L}$ coincides with
$E^{\emptyset ,\delta}_{KP}$ for finite $\delta$ and
$Aut_{L}(M)$ is closed in $Aut(M)$.

\bigskip

If one is willing to allow many-sorted $\omega$-categorical
structures, then the non-$G$-compact theory obtained in
\cite{clpz} (denoted by $(M_i )_{i\in\omega}$) is already
$\omega$-categorical.
If one insists on a one-sorted $\omega$-categorical structure,
one should build a structure $N$ with maps $f_n$ from
${N \choose n}$ (the set of $n$-element subsets of $N$) to $M_n$
in a suitably generic fashion so that $N$ is
$\omega$-categorical.
This is roughly what we do below.
It turns out that this idea is realized by modifying
of some construction (Example 3.4) from \cite{ivma}.
\parskip0pt

Moreover we have found that the final structure does not have
AZ-enumerations, therefore it is the first example of an
$\omega$-categorical structure without AZ-enumerations
(answering a question from \cite{AlCh}).
The latter notion was introduced by Hrushovski in \cite{hrus}
(in different terminology) as a technical tool allowing him
to solve several basic problems concerning $\omega$-stable
$\omega$-categorical structures.
This notion refines the notion of {\em nice emumerations}
introduced earlier by Ahlbrandt and Ziegler in \cite{az} and
since then involved in many questions of model theory.
\parskip0pt

Analyzing our construction we have found that the absence
of an AZ-enumeration can be already obtained for some reducts
of our structure.
This is the reason why the next section of the paper is devoted
to the easiest example of an $\omega$-categorical structure
without AZ-enumerations we can build.
The $\omega$-categorical non-$G$-compact structure which we have
found, is an expansion of this example and will be described in
Section 2.

\section{Enumerations}
In this section we give some easy construction of an
$\aleph_{0}$-categorical structure which does not have
AZ-enumerations.
This is a simplified version of the main construction of
the paper.
It can be considered as a nice warm up.

A linear ordering $\prec$ of a countable structure $M$ is
called an {\em AZ-enumeration} of $M$ if it has order-type
$\omega$ and for any $n\ge 1$ it satisfies the following
property:
\begin{quote}
whenever $\bar{b}_i$, $i<\omega$, is a sequence of
$n$-tuples from $M$, there exist some $i<j<\omega$ and
a $\prec$-preserving elementary map $f:M\rightarrow M$
such that $f(\bar{b}_i )=\bar{b}_j$.
\end{quote}

Let $L_0 =\{ E_n : 2< n<\omega \}$ be a first-order language,
where each $E_n$ is a relational symbol of arity $2n$.
Let $K$ be the class of all finite $L_0$-structures $C$ where
each relation $E_n (\bar{x},\bar{y})$ determines an equivalence
relation on the set (denoted by ${C \choose n}$) of unordered
$n$-element subsets of $C$.
In particular we have that $K$ satisfies the sentence
$$
\forall \bar{x}\bar{y} (E_n (x_1 ,...,x_n ,y_1 ,...,y_n )\rightarrow
\bigwedge \{ E_{n}(y_1 ,...,y_n ,x_{\sigma (1)},...,x_{\sigma (n)}):
\sigma \in Sym (n)\}).
$$
For $C\in K$ and $n>|C|$ we put that no $2n$-tuple from $C$
satisfies $E_n (\bar{x},\bar{y})$.
It is easy to see that $K$ is closed under taking
substructures and the number of isomorphism types of
$K$-structures of any finite size is finite.
\parskip0pt

To verify {\em the amalgamation property for} $K$,
given $A,B_{1},B_{2} \in K$ with $B_{1}\cap B_{2}=A$,
define $C\in K$ as $B_{1}\cup B_{2}$.
The relations $E_{n}$, $n\le |B_{1}\cup B_{2}|$, can be
easily defined so that $C\in K$ and $B_i <C$.
To be more precise we will obey the following rules.
When $n\le |B_{1}\cup B_{2}|$ and
$\bar{a}\in {B_1 \choose n} \cup {B_2 \choose n}$ we put
that the $E_n$-class of $\bar{a}$ in $C$ is contained in
${B_1 \choose n} \cup {B_2 \choose n}$.
We also assume that all $n$-tuples meeting both
$B_{1}\setminus B_{2}$ and $B_{2}\setminus B_{1}$
are pairwise equivalent with respect to $E_{n}$.
In particular if $n\ge max(|B_{1}|,|B_{2}|)$ we put that all
$n$-tuples from $C$ are pairwise $E_{n}$-equivalent.
\parskip0pt

It is easy to see that this amalgamation also works for
the joint embedding property.\parskip0pt

Let $M$ be the countable universal homogeneous structure for $K$.
It is clear that in $M$ each $E_n$ defines infinitely many
classes and each $E_n$-class is infinite.

\begin{prop}
The structure $M$ is $\aleph_{0}$-categorical and does not have
any AZ-enumeration.
\end{prop}

{\em Proof.}
Since for each $n$ the number of finite structures of $K$ of
size $n$ is finite, the structure $M$ is
$\aleph_{0}$-categorical and admits elimination of quantifiers
(by Fraiss\'{e}'s theorem).
Now for a contradiction suppose that there is an ordering $\prec$
defining an AZ-enumeration of $M$.
We will define an infinite sequence of triples
$a_n \prec b_n \prec c_n$, $n\in \omega$, satisfying the
following conditions.
Let $a_n$ be always the minimal $\prec$-element of $M$.
For $n>1$ the elements $b_n$ and $c_n$ are chosen so that for any
$(n-1)$-tuple of the form $x_1 \prec x_2 \prec ...\prec x_{n-1}$
with $x_{n-1}\prec b_n$ the tuple $(x_1 ,x_2 ,...,x_{n-1},c_n )$
is $E_n$-equivalent with some $n$-tuple $\bar{y}$ satisfying
$y_1 \prec ...\prec y_n \prec b_n$.
On the other hand we also demand that for each $j<n$, any $j$-tuple 
of the form $d_1 \prec d_2 \prec ...\prec d_{j-1}\prec c_n$ 
with $d_{j-1}\prec b_n$ is not $E_j$-equivalent with any $j$-tuple 
$y_1 \prec ...\prec y_j$ with $y_j \prec b_n$. \parskip0pt

The triples $(a_n ,b_n ,c_n )$ can be defined by induction.
Let $a_1 \prec b_1 \prec c_1$ be the initial 3-element
$\prec$-segment of $M$.
At step $n$ we just take $b_n$ as the next element enumerated
after $c_{n-1}$.
To define $c_n$ consider the substructure of $M$ defined on
$D=\{ x : x\prec b_n \}$.
We embed $D$ into some $K$-structure $D\cup \{ c\}$
such that for each $j<n$ all tuples $(y_1 ,...,y_{j-1},c)$ with 
$y_1 \prec y_2 \prec...\prec y_{j-1}\prec b_n$
form an $E_j$-class which does not meet any $j$-tuple from $D$.
We also demand that each $n$-tuple of $D\cup \{ c\}$ is 
$E_n$-equivalent with an $n$-tuple of $D$.
Since $M$ is universal homogeneous, the element $c$ can be
found in $M$.
Let $c_n$ be the element of $M$ with $D\cup \{ c_n \}$
isomorphic with $D\cup\{ c\}$ over $D$ and having the
minimal number with respect to $\prec$. \parskip0pt

If $f:M\rightarrow M$ is a $\prec$-preserving elementary
map taking $(a_i ,b_i ,c_i )$ to $(a_j ,b_j ,c_j )$,
then by the definition of $c_j$ any $i$-tuple
of $\{ x:x\prec b_j \}\cup\{ c_j\}$ with $c_j$ is not
$E_i$-equivalent with any tuple of $\{ x: x\prec b_j \}$.
By the definition of $c_i$ this is impossible.
Therefore we have a contradiction with the definition
of an AZ-enumeration.

\section{The main example}
We build our structure by a generalized Fraiss\'{e}'s
construction, appealing to Theorem 2.10 of \cite{evans}, p.44.
We now recall that material.\parskip0pt

Let $L$ be a relational language and let $\mathcal{C}$ be
a class of finite $L$-structures.
Let $\mathcal{E}$ be a class of embeddings
$\alpha :A\rightarrow B$ (where $A,B\in \mathcal{C}$) such that
any isomorphism $\delta$ between $\mathcal{C}$-structures
(from $Dom(\delta )$ onto $Range(\delta )$) is in $\mathcal{E}$,
the class $\mathcal{E}$ is closed under composition and
the following property holds:
\begin{quote}
if $\alpha :A\rightarrow B$ is in $\mathcal{E}$ and
$C\subseteq B$ is a substructure in $\mathcal{C}$ such that
$\alpha (A)\subseteq C$, then the map obtained by restricting
the range of $\alpha$ to $C$ is also in $\mathcal{E}$.
\end{quote}
We say that a structure $A\in \mathcal{C}$ is a {\em strong}
substructure of an $L$-structure $M$ if $A\subseteq M$ and
any inclusion $A\subseteq B$ with $B\in \mathcal{C}$ and
$B\subseteq M$ is an $\mathcal{E}$-embedding.
We call an embeddings $\rho :C\rightarrow M$ {\em strong} if
$C\in\mathcal{C}$ and $\rho (C)$ is a strong substructure of $M$.
\parskip0pt

Theorem 2.10 of \cite{evans} states that if

\begin{quote}
(a) the number of isomorphism types of $\mathcal{C}$-structures
of any finite size is finite; \parskip0pt

(b) the class $\mathcal{E}$ satisfies the joint embedding
property and the amalgamation property and \parskip0pt

(c) there is a function $\theta$ on the natural numbers such
that any $L$-structure $C$ embeds into some $A\in \mathcal{C}$
of size $\le\theta (|C|)$ such that any embedding from $A$ to
a $\mathcal{C}$-structure is strong;
\end{quote}
then there exists a countably categorical $L$-structure $M$
such that $M$ is {\em generic}, i.e.

\begin{quote}
(a') $\mathcal{C}$ is the class of all strong substructures
of $M$; \parskip0pt

(b') $M$ is a union of a chain of $\mathcal{E}$-embeddings and
\parskip0pt

(c') if $A$ is a strong substructure of $M$ and
$\alpha :A\rightarrow B$ is in $\mathcal{E}$
then $B$ is strongly embeddable into $M$ over $A$.
\end{quote}
Moreover any isomorphism between strong finite substructures
of $M$ extends to an automorphism of $M$.\parskip0pt

Let $L=\{ E_{n},K_{n},R_{n}: 2<n\in\omega \}$ be a first-order
language, where each $E_{n}$ and $R_{n}$ is a relational symbol
of arity $2n$ and each $K_{n}$ has arity $3n$.
The structure $M$ which is anounced in Introduction, will be
built by the version of Fraiss\'{e}'s construction described
above.
We first specify a class $K$ of finite $L$-structures, which
will become the class of all finite $L$-substructures of $M$.
\parskip0pt

As in Section 1 in each $C\in K$ each relation $E_{n}$
determines an equivalence relation on the set (denoted by
${C \choose n}$) of unordered $n$-element subsets of $C$.
The relations $R_{n}$ are irreflexive.
The $R_{n}$-arrows respect $E_{n}$,
$$
\forall \bar{x},\bar{y},\bar{u},\bar{w} (E_{n}(\bar{x},\bar{y})
\wedge E_{n}(\bar{u},\bar{w})\wedge R_{n}(\bar{x},\bar{u})
\rightarrow R_{n}(\bar{y},\bar{w})),
$$
and define a partial 1-1-function on ${C \choose n} /E_{n}$.
\parskip0pt

Every $K_{n}$ is interpreted by a circular order
\footnote{a twisted around total order with the natural ternary
relation induced by $<$}
on the set of $E_{n}$-classes.
Therefore we take the axiom
$$
\forall \bar{x},\bar{y},\bar{z},\bar{u},\bar{v},\bar{w}
(E_{n}(\bar{x},\bar{y})\wedge E_{n}(\bar{u},\bar{w})\wedge
E_{n}(\bar{z},\bar{v})\wedge K_{n}(\bar{x},\bar{z},\bar{u})
\rightarrow K_{n}(\bar{y},\bar{v},\bar{w})).
$$
and the corresponding axioms of circular orders.
We also take some axioms connecting $K_{n}$ and $R_n$:
$$
R_{n}(\bar{x},\bar{y}) \wedge R_{n}(\bar{y},\bar{z}) \rightarrow
K_{n}(\bar{x},\bar{y},\bar{z});
$$
$$
\forall \bar{v}_{1},\bar{v}_{2},\bar{v}_{3},\bar{w}_{1},
\bar{w}_{2},\bar{w}_{3}
(\bigwedge_{i\le 3} R_{n}(\bar{v}_{i},\bar{w}_{i})\rightarrow
(K_{n}(\bar{v}_{1},\bar{v}_{2},\bar{v}_{3})\leftrightarrow
K_{n}(\bar{w}_{1},\bar{w}_{2},\bar{w}_{3}))).
$$
These axioms say that $R_n$ defines a partial automorphism of
the circular order induced by $K_n$ on ${C \choose n}/E_n$.
Our final axioms state that this partial automorphism
admits an extension to a 1-1-function $f$ (on some larger
domain) such that $f^{n}$ is identity on its domain, but for
each $V\in {C \choose n} /E_{n}$ and $m\not=0$ with $m<n$ we
have $f^{m}(V)\not= V$.
These conditions can be written by an infinite set of
universal first-order formulas (which forbid all inconsistent
situations). \parskip0pt

It is easy to see that the class $K$ is closed under taking
substructures.
It is noted in \cite{ivma} that the class of reducts of
$K$-structures to $\{ E_{n},K_{n}:n>2\}$ has the amalgamation
property.
Then the example given in \cite{ivma} is just the universal
homogeneous structure defined by these reducts.
It is shown in \cite{ivma} that it does not admit strongly
determined types over any finite set. \parskip0pt

On the other hand $K$ does not satisfy the amalgamation property.
We now describe a cofinal subclass $\mathcal{C}\subset K$ with
the amalgamation property.
The variant of Fraiss\'{e}'s theorem described in the
beginning of the section will be applied to this subclass
$\mathcal{C}$. \parskip0pt

We say that a structure $A\in K$ is {\em strong}, if for every
$n\in\omega$ all elements of ${A \choose n}$ are pairwise
equivalent with respect to $E_{n}$ or for any
$\bar{a}\in {A \choose n}$
there is a sequence $\bar{a}_{1}(=\bar{a}),...,\bar{a}_{n}$ of
pairwise non-$E_{n}$-equivalent tuples from ${A \choose n}$
such that $(\bar{a}_{i},\bar{a}_{i+1})\in R_{n}$,
$1\le i\le n-1$, and $(\bar{a}_{n},\bar{a}_{1})\in R_{n}$.
Let $\mathcal{C}$ be the class of all strong structures from $K$.
\parskip0pt

Let us show that $\mathcal{C}$ {\em is cofinal in} $K$.
Let $C\in K$ and $|{C \choose n} /E_{n}|\not= 1$.
For a $\bar{c}\in {C \choose n}$ witnessing that $C$
is not strong (in particular there is no sequence
$\bar{c}_{1}(=\bar{c}),...,\bar{c}_{n}$ of tuples from
${C \choose n}$ such that
$(\bar{c}_{i},\bar{c}_{i+1})\in R_{n}$, $1\le i\le n-1$,
and $(\bar{c}_{n},\bar{c}_{1})\in R_{n}$) let $C\cup C'\in K$
be a structure defined on the disjoint union of $C$ and $C'$,
which contains $C$ as a substructure and has the property
that there is a sequence $\bar{c}_{2},...,\bar{c}_{n}$
of pairwise non-$E_{n}$-equivalent tuples from
${ C\cup C' \choose n}$ such that $(\bar{c}_{i},\bar{c}_{i+1})\in R_{n}$,
$1\le n-1$, and $(\bar{c}_{n},\bar{c}_{1})\in R_{n}$
(it can happen that some $\bar{c}_{i}$ are in $C$).
We also assume that for each $m$ all $m$-tuples from
${C\cup C' \choose m} \setminus \{\bar{c}_{2},...,\bar{c}_{n}\}$
meeting $C'$ are $E_{m}$-equivalent and moreover they are
$E_{m}$-equivalent to some fixed $\bar{c}\subset C$ if
$|C|\ge m$.
Note that the number of tuples witnessing that $C\cup C'$
is not strong is less than that for $C$.
At the second step we repeat this construction for the next
tuple in $C\cup C'$.
As a result we obtain some $C\cup C'\cup C''$.
Continuing this procedure we obtain in finitely many steps
a strong structure $C\cup C'\cup ...\cup C^{(k)}$.
\parskip0pt

It is worth noting that at every step we can arrange
that $|C^{(l)}|\le |C|^{2}$.
On the other hand the number of steps is not greater
than $2^{|C|}$.
As a result we see that the size of the structure obtained does
not exceed  $2^{|C|}|C|^{2}$. \parskip0pt

We now verify {\em the amalgamation (and the joint embedding)
property for} $\mathcal{C}$.
Given $A,B_{1},B_{2} \in \mathcal{C}$ with $B_{1}\cap B_{2}=A$,
define $C\in \mathcal{C}$ as $B_{1}\cup B_{2}$.
The relations $E_{n},R_{n},K_{n}$,
$n\le |B_{1}\cup B_{2}|$, are defined so that $C\in K$
and the following conditions hold.
Let $n\le |B_{1}\cup B_{2}|$.
We put that all $n$-tuples meeting both $B_{1}\setminus B_{2}$
and $B_{2}\setminus B_{1}$ are pairwise equivalent with respect
to $E_{n}$.
We additionally demand that they are equivalent to some tuple
from some $B_{i}$, $i\in \{ 1,2\}$, if $n\le max(|B_{1}|,|B_{2}|)$.
If for some $i\in \{ 1,2\}$, $|{ B_{i} \choose n} /E_{n}|=1$,
then we put that all $n$-tuples $\bar{c}\in B_{1}\cup B_{2}$
meeting $B_{i}$ are pairwise $E_{n}$-equivalent.
We additionally assume that they are equivalent
to some tuple from $B_{3-i}$ if $n\le |B_{3-i}|$.
If
$|{B_{1}\choose n} /E_{n}|\not= 1\not=|{B_{2}\choose n} /E_{n}|$,
$\bar{c}_{1},\bar{c}_{2},...,\bar{c}_{n}$ is an
$R_{n}$-cycle in $B_{i}$ and
$\bar{c}'_{1},\bar{c}'_{2},...,\bar{c}'_{n}$ is an $R_{n}$-cycle
in $B_{3-i}$ with $(\bar{c}_{1},\bar{c}'_{1})\in E_{n}|_{A}$,
then we define $E_{n}$ and $R_{n}$ so that for all $i\le n$,
$(\bar{c}_{i},\bar{c}'_{i})\in E_{n}$
(since $A\in \mathcal{C}$ this can be easily arranged).
If $n\ge max(|B_{1}|,|B_{2}|)$ then all $n$-tuples from
$C$ are pairwise $E_{n}$-equivalent.
We assume that $E_n$ is the minimal equivalence relation
satisfying the conditions above. \parskip0pt

We can now define the circular orderings $K_{n}$ on $C$.
There is nothing to do if $|{C \choose n} /E_{n}|=1$.
In the case when for some $i=1,2$,
$|{B_{i} \choose n} /E_{n}| = 1$, the relation $K_{n}$ is
defined by its restriction to $B_{3-i}$.
When $|{B_{1} \choose n} /E_{n}|\not= 1\not= |{B_{2} \choose n} /E_{n}|$,
the ordering $K_{n}$ on an $R_{n}$-cycle corresponds
to the relation $R_{n}$.
Thus for any $R_{n}$-cycle in ${C \choose n} /E_{n}$ having
representatives both in $B_{1}$ and in $B_{2}$ (see the previous
paragraph) the definition of $K_{n}$ does not depend on
the choice of representatives.
In the case when such a cycle exists we fix an element
$V\in {C \choose n} /E_{n}$ of this cycle and $V'$ with
$(V,V')\in R_{n}$.
Then amalgamate the linear orderings between $V$ and $V'$
in $(B_{1},K_{n})$ and in $(B_{2},K_{n})$
(over the set of $E_n$-classes having representatives both in
${B_1 \choose n}$ and ${B_2 \choose n}$).
This defines $K_{n}$ on ${C \choose n} /E_{n}$. \parskip0pt

If all $R_{n}$-cycles in ${C \choose n} /E_{n}$ having
representatives in $B_{1}$ do not have representatives in
$B_{2}$ we fix elements $V_{1},V_{2}\in {C \choose n} /E_{n}$
which represent $R_{n}$-cycles in $B_{1}$ and in $B_{2}$
respectively.
If $(V_{1},V'_{1})\in R_{n}$ and $(V_{2},V'_{2})\in R_{n}$
we put (amalgamating the corresponding linear orderings) that
all elements between $V_{2}$ and $V'_{2}$ (including $V_{2}$
and not including $V'_2$) are between $V_{1}$ and $V'_{1}$
(with the same direction).
This defines $K_{n}$ on ${C \choose n} /E_{n}$. \parskip0pt

Let $\mathcal{E}$ be the class of all embeddings between
(strong) structures from $\mathcal{C}$.
It is clear that the number of isomorphism types of
$L$-structures of any finite size is finite (thus condition
(a) above is satisfied).
We have already noticed that
(b) strong embeddings of $\mathcal{C}$-structures satisfy the
joint embedding property and the amalgamation property.
We have also shown that
\begin{quote}
(c) the function $\theta (n)=2^{n}\cdot n^2$ satisfies
the property that any $L$-structure $C$ embeds into some
$A\in \mathcal{C}$ of size $\le\theta (|C|)$ (and any embedding
from $A$ to a $\mathcal{C}$-structure is strong).
\end{quote}
By the version of Fraiss\'{e}'s Theorem from \cite{evans}
described in the beginning of the section, there exists
a countably categorical structure $M$ such that $M$ is
{\em generic}: (a') $\mathcal{C}$ is the
class of all strong substructures of $M$, (b') $M$ is a union of
a chain of strong embeddings and (c') if $A$ is a strong
substructure of $M$ and $\alpha :A\rightarrow B$ is a strong
embedding with $B\in \mathcal{C}$, then $B$ is embeddable into
$M$ over $A$.
Moreover any isomorphism between strong finite substructures of
$M$ extends to an automorphism of $M$.\parskip0pt

We now want to prove that the theory $Th(M)$ is not G-compact.
To simplify notation below we often replace formulas
of the form $K_{n}(\bar{a},\bar{b},\bar{c})$ by expressions
$\bar{a}<_{n}\bar{b}<_{n}\bar{c}$.
We also apply $\le_{n}$ to $E_n$-classes when it is convenient
to identify an $E_{n}$-class with its representative neglecting
the difference.
By $f(\bar{c})$ we denote some (any) $\bar{c}'$ with
$R_{n}(\bar{c},\bar{c}')$ (we do not write $f_{n}(\bar{c})$
because $n$ equals the length of $\bar{c}$).\parskip0pt

The following lemma is a standard application of genericity.

\begin{lem}
Let $M$ be a generic structure for $\mathcal{C}$ and
$M\models R_{n}(\bar{a},\bar{b})$.
Then the linear ordering induced by $K_{n}$ on the set of
$E_{n}$-classes of $\{\bar{c}: K_{n}(\bar{a},\bar{c},\bar{b})\}$
is dense and without endpoints.
\end{lem}

We now describe our main tool for non-G-compactness.

\begin{lem}\label{fir}
Let $M$ be a generic structure for $\mathcal{C}$.
Let $\bar{c}_{1}$ and $\bar{c}_{2}\in M$ enumerate strong
substructures of the same type over $\emptyset$ such that
\parskip0pt

(a)
$M\models (\bar{c}_{1}\le_{n}\bar{c}_{2}<_{n}f(\bar{c}_1 ))$;
\parskip0pt

(b) the tuple $\bar{c}_{1}\bar{c}_{2}$ enumerates a strong
substructure $D$ where for every $m$, every $R_{m}$-cycle
in ${D \choose m} /E_{m}$ is already realized by
$E_{m}$-equivalence in either ${\bar{c}_{1} \choose m} /E_{m}$
or ${\bar{c}_{2} \choose m} /E_{m}$
and can not be realized in $\bar{c}_1$ and $\bar{c}_2$
simultaneously;\parskip0pt

(c) for any pair of subtuples $\bar{c}'_{1}\subseteq \bar{c}_{1}$
and $\bar{c}'_{2}\subseteq \bar{c}_{2}$
(say $|\bar{c}'_1 |=|\bar{c}'_2 |=m$) representing the same
places in $\bar{c}_{1}$ and $\bar{c}_{2}$, we have
$$
M\models (\bar{c}'_{1}\le_{m}\bar{c}'_{2}<_{m}f(\bar{c}'_{1}))
\vee (\bar{c}'_{2}\le_{m}\bar{c}'_{1}<_{m}f(\bar{c}'_{2})).
$$
Then there is an elementary substructure $N$ of $M$ such that
$\bar{c}_{1}\bar{c}_{2}\cap N=\emptyset$ and
$tp(\bar{c}_{1}/N)= tp(\bar{c}_{2}/N)$.
\end{lem}

{\em Proof.} Consider a chain of strong
embeddings $C_{1}<C_{2}< ...$ such that $M=\bigcup C_{i}$.
We build a chain $C'_{1}<C'_{2}< ...$ of
strong substructures of $M$ together with an increasing chain
$\varepsilon_{1} \subseteq \varepsilon_2 \subseteq ...$ of
isomorphisms $\varepsilon_i :C'_{i}\rightarrow C_i$
such that for every $i$, $D\cap C'_{i}=\emptyset$ and
$tp(\bar{c}_{1}/C'_{i})= tp(\bar{c}_{2}/C'_{i})$.
If such a chain exists then $N=\bigcup C'_{i}$ is also
a generic structure and by Tarski-Vaught test $N$ is
an elementary substructure of $M$ (in fact $Th(M)$ is
model complete).\parskip0pt

The condition
$tp(\bar{c}_{1}/C'_{i})= tp(\bar{c}_{2}/C'_{i})$
will be satisfied as follows.
At every step of our construction we find $C'_{i}$ so that
$DC'_{i}$ is a strong substructure of $M$, where all
$E_m$-classes from $R_m$-cycles of $D$ remain the same as in $D$.
Having this we can additionally arrange that $\bar{c}_1 C'_{i}$
and $\bar{c}_2 C'_{i}$ are strong substructures of $DC'_{i}$
which are isomorphic over $C'_{i}$ with respect to the map
$\bar{c}_1 \rightarrow \bar{c}_2$.
As any isomorphism of strong substructures extends to an
automorphism of $M$ we will see that
$tp(\bar{c}_{1} /C'_{i})=tp(\bar{c}_{2} /C'_{i})$. \parskip0pt

We may assume that for all $k$,
$|\bar{c}_{1}\bar{c}_{2}|+|C_{k}|<|C_{k+1}|$ and for every
$m$ with an $R_{m}$-cycle in the structure $D$ there is an
$R_{m}$-cycle in $C_{1}$ (not necessarily $E_m$-equivalent
with the former one).
The latter assumption will guarantee that in the construction
below $E_{m}$-classes from $R_{m}$-cycles of
$D$ remain the same in the extended structures
$DC'_{k}$. \parskip0pt

The existence of a chain $C'_{1}<C'_{2}< ...$ as above
will be shown by induction.
Assume that there are strong substructures
$C'_{1}<...<C'_{k}<M$ such that for every $j\le k$,
$\bar{c}_{1}C'_{j}$, $\bar{c}_{2}C'_{j}$ and
$DC'_{j}$ are strong substructures and
$tp(\bar{c}_{1}/C'_{j})= tp(\bar{c}_{2}/C'_{j})$.
Define $B\in \mathcal{C}$ (which will be a copy of $DC'_{k+1}$)
as the quotient of the disjoint union $C_{k+1}\cup DC'_{k}$
by the isomorphism
$\varepsilon_{k}: C'_{k}\rightarrow C_{k}<C_{k+1}$
(we identify images with their preimages).
The relations $E_{m}$ and $R_{m}$, $m\le |B|$, are
defined as in the amalgamation procedure described
above (we replace $B_{1}$ by $C_{k+1}$, $B_{2}$
by $DC'_{k}$ and $A$ by $C'_{k}$).
We now make a small modification in this procedure: we put that
all $m$-tuples meeting both $C_{k+1}\setminus C_{k}$ and $D$
are pairwise equivalent with respect to $E_{m}$ and we
additionally demand that they are equivalent
to some tuple from $C_{k+1}$ if $m\le |C_{k+1}|$
(we always assume that $|DC'_{k}|<|C_{k+1}|$).
If $|{C_{k+1} \choose m} /E_{m}|=1$, then we put that all
$m$-tuples meeting $C_{k+1}$ are pairwise
$E_{m}$-equivalent and if $m\le |DC'_{k}|$
we put that these tuples are equivalent to some
(any) $m$-tuple from $DC'_{k}$.
It is worth noting here that at Step 1 (where we assume that
$C_{0}=\emptyset$) the condition $|{C_{1} \choose m} /E_{m}|=1$
implies $|{D \choose m} /E_{m}|=1$ and then we define
$|{DC_{1} \choose m} /E_{m}|=1$.
Since $|{C_{1} \choose m} /E_{m}|=1$ follows from
$|{C_{k+1} \choose m} /E_{m}|=1$,
we have that $|{C_{k+1} \choose m} /E_{m}|=1$
always implies $|{B \choose m} /E_{m}|=1$.\parskip0pt

It is now easy to see that the amalgamation procedure
in the form above guarantees that:

(i) all $E_{m}$-classes from $R_{m}$-cycles of
$D$ remain the same in the extended structure $B$;
\parskip0pt

(ii) if $D$ contains $R_m$-cycles, there is a unique $E_m$-class
from $B$ which does not have any element which is a subtuple of
some $\bar{c}_i$, $i\in\{ 1,2\}$, but contains some tuples which
intersect $D$; this is the class containing all possible tuples
from $B$ meeting both $C_{k+1}\setminus C_k$ and $D$; \parskip0pt

(iii) $\bar{c}_{1}$ and $\bar{c}_{2}$ realize in $B$ the same
quantifier-free type over $C_{k+1}$ with respect to
the sublanguage $\{ E_{n},R_{n}: 2<n\in\omega\}$.\parskip0pt

We also modify the construction of the circular ordering $K_{m}$
on $B$.
There is nothing to do if $|{DC'_k \choose m} /E_{m}|=1$ (then
$K_m$ on $B$ is determined by its restriction to $C_{k+1}$).
Assume $|{DC'_k \choose m} /E_{m}|\not= 1$.
At the first step (when we amalgamate $D$ with $C_1$)
find $V_{1},V_{2}\in {B \choose m} /E_{m}$ which represent
some $R_{m}$-cycles in $C_{1}$ and in $D$ respectively.
If $(V_{1},V'_{1})\in R_{m}$ and
$(V_{2},V'_{2})\in R_{m}$ we put that all elements between
$V_{2}$ and $V'_{2}$ (including $V_{2}$ and not including $V'_2$)
are between $V_{1}$ and $V'_{1}$ (with the same direction).
After appropriate amalgamation we obtain $K_{m}$ on
${B \choose m} /E_{m}$.
To guarantee that $\bar{c}_{1}$ and $\bar{c}_{2}$ have
the same type over $C_{1}$ in $B$ we must only consider
$E_{m}$-classes representing $R_{m}$-cycles with
tuples from $\bar{c}_{i}$, $i\in \{ 1,2\}$ (by (ii)).
Then amalgamating the linear orderings $(V_{1},V'_{1})$
and $(V_{2},V'_{2})$ as above we put that $E_{m}$-classes
corresponding to matched subtuples of $\bar{c}_{i}$,
$i\in \{ 1,2\}$ (these classes are the same as
the corresponding ones in $D$),  are not separated by
any $E_{m}$-class having representatives meeting $C_1$.
This can be done by the last assumption of the lemma.
\parskip0pt

At later steps note that if $C'_k$ does not have $R_m$-cycles,
then there are no $R_m$-cycles in $D$.
Then our definition of $E_m$ implies
$|{DC'_k \choose m} /E_{m}|= 1$, a contradiction.
We see that there are $R_{m}$-cycles in ${B \choose m} /E_{m}$
having representatives (by $E_m$) both in $C_{k+1}$ and in
$DC'_{k}$: for example these are $R_m$-cycles occurring in
$C'_k$.\parskip0pt

If an $R_{m}$-cycle in ${B \choose m} /E_{m}$ has
representatives both in $C_{k+1}$ and in $DC'_{k}$, then
the definition of $K_{m}$ on this cycle does not depend on
the choice of representatives.
In this case we fix an element $V\in {B \choose m})/E_{m}$
of such a cycle and amalgamate the linear orderings between $V$
and $V'$ with $(V,V')\in R_{m}$ in $C_{k+1}$ and in $D C'_{k}$
(as in the process of amalgamation described above).
We again put that in these linear orderings
$E_{m}$-classes corresponding to matched subtuples of
$\bar{c}_{i}$, $i\in \{ 1,2\}$, are not separated by any
$E_{m}$-class meeting $B\setminus D$.
Here we again apply the last assumption of the lemma
and inductive hypotheses (in particular
$tp(\bar{c}_{1}/C'_{k})=tp(\bar{c}_{2}/C'_{k})$).
As a result we obtain that $\bar{c}_{1}$ and
$\bar{c}_{2}$ have the same quantifier-free type over
$C_{k+1}$ in $B$. \parskip0pt

Using the fact that $M$ is generic we embed $B$ into $M$
over $D C'_{k}$.
The image of $C_{k+1}$ is the required structure $C'_{k+1}$
and $\varepsilon_{k+1}$ is the converse map.
It is clear that $D C'_{k+1}$ is strong.
To prove that
$tp(\bar{c}_{1}/C'_{k+1})= tp(\bar{c}_{2}/C'_{k+1})$ in $M$
it suffices to show that structures $\bar{c}_{1}C'_{k+1}$
and $\bar{c}_{2}C'_{k+1}$ are strong and isomorphic over
$C'_{k+1}$.\parskip0pt

We have already noticed that the condition
$|{C_{i} \choose m} /E_{m}|=1$ where $i\le k+1$, implies that
$|{C_{1} \choose m} /E_{m}|=1$, $|{D \choose m} /E_{m}|=1$ and
that for all $j\le i$ there is a unique $E_{m}$-class over
$DC'_{j}$.
As a result if $|{\bar{c}_{1} \choose m} /E_{m}|\not= 1$
(which is equivalent to $|{D \choose m} /E_{m}|\not= 1$), then
the $E_{m}$-classes of $R_{m}$-cycles from ${D \choose m}$
remain the same in $D C'_{k+1}$.
Thus for such $m$ any $E_{m}$-class $V$ in
$\bar{c}_{1}C'_{k+1}$ belongs to an $R_{m}$-cycle
(if $V$ contains a tuple meeting $C'_{k+1}$ then $V$ belongs
to an $R_{m}$-cycle defined in $C'_{k+1}$).
If $|{\bar{c}_{1} \choose m} /E_{m}|= 1$, then the
$R_{m}$-cycles in $\bar{c}_{1}C'_{k+1}$ are defined by
$R_{m}$-cycles over $C'_{k+1}$.
This shows that $\bar{c}_{1}C'_{k+1}$ is strong.
Similar arguments and the definition of $K_{m}$ imply that
$\bar{c}_{1}C'_{k+1}$ and $\bar{c}_{2}C'_{k+1}$ are isomorphic
over $C'_{k+1}$.
$\square$

\bigskip

We also need the following lemma.

\begin{lem} \label{sec}
Let tuples $\bar{a}$ and $\bar{b}$ from $M$ enumerate
strong substructures of the same type over $\emptyset$.
If they have the same type over some elementary substructure
of $M$ then for any pair of subtuples $\bar{a}'$
and $\bar{b}'$ (of length $l$) representing the same places in
$\bar{a}$ and $\bar{b}$, the following condition holds:
$$
(\bar{a}'\le_{l}\bar{b}'<_{l}f(\bar{a}'))
\vee(\bar{b}'\le_{l}\bar{a}'<_{l}f(\bar{b}')).
$$
\end{lem}

{\em Proof.} Assume that the condition does not hold for
subtuples $\bar{a}'$ and $\bar{b}'$.
For any elementary substructure $N$ there are tuples $\bar{c}$
and $\bar{d}\in N$ which are $K_{l}$-between $\bar{a}'$ and
$f(\bar{a}')$ and $\bar{b}'$ and $f(\bar{b}')$ respectively.
Then $\bar{a}'$ and $\bar{b}'$ do not have the same
type over $\bar{c}\bar{d}$. $\square$

\begin{thm} \label{thr}
The theory $Th(M)$ is not G-compact.
\end{thm}

{\em Proof.}
We present the structure $M$ as two sequences of
strong embeddings $C_{1}<C_{2}<...$ and $C'_{1}<C'_{2}<...$ such
that for every  $i$,
$tp(C_{1},...,C_{i}/\emptyset )= tp(C'_{1},...,C'_{i}/\emptyset )$
(under appropriate enumerations of $C_{j},C'_{k}$) and there are
$D_{-i},...,D_{-1},D_{1},...,D_{i}\subseteq M$ such that all
pairs
$$
(D_{-i},D_{-i+1}),...,(D_{-2},D_{-1}),(D_{-1},C_{4i}),
(C_{4i},D_{1}),...,(D_{i-1},D_{i})
$$
belong to $R_{|C_{4i}|}$ and the corresponding
$K_{|C_{4i}|}$-intervals do not contain $C'_{4i}$.
The first sequence can be chosen arbitrary.
Then $|C_{4i}|\ge 4i$ and $C_{4i}$ belongs to a
$R_{|C_{4i}|}$-cycle of length $\ge 4i$.
The existence of the second sequence can be obtained by
induction where at every step we apply genericity of $M$.
\parskip0pt

Let $\alpha$ be an automorphism of $M$ taking every $C_{i}$ to
$C'_{i}$.
If $\alpha$ is a product of $n$ automorphisms fixing
elementary substructures of $M$, then by Lemma \ref{sec} we have
a contradiction with the existence $D_{i}$ as above for $C_{4n}$
and $C'_{4n}$.
As a result we have that $\alpha$ is not a Lascar
strong automorphism. \parskip0pt

On the other hand for any finite map of the form
$\alpha |_{C_{i}}$ find $\bar{c}_{1}$ (enumerating $C_{i}$),
$\bar{c}_{2},...,\bar{c}_{k}$ ($\bar{c}_{k}$ enumerates $C'_{i}$
in the appropriate way) such that any pair
$\bar{c}_{j},\bar{c}_{j+1}$ satisfies the conditions of
Lemma \ref{fir}.
To find such tuples we apply genericity of $M$: we can ensure
that tuples of each pair do not have common elements;
then it is easy to arrange (by amalgamation) that every
pair forms a strong structure as in Lemma \ref{fir}.
Applying the lemma to these pairs we obtain that
$\alpha |_{C_{i}}$ can be
presented as a restriction of a Lascar strong automorphisms
(generated by automorphisms fixing elementary substructures).
We see that $\alpha$ belongs to the closure of the group of
Lascar strong automorphisms in $Aut(M)$.
This implies that $Th(M)$ is not G-compact. $\square$

\bigskip

The example of the previous section is a reduct of the
structure obtained above.
Since any AZ-enumeration of a structure $M$ is
an AZ-enumeration of any its reduct, we see that $M$
does not have AZ-enumerations. \parskip0pt

We finish the paper by a remark concerning diameters of
Lascar strong types.
They are defined in \cite{clpz} as follows.
For $\bar{a}$ and $\bar{b}$ let $d(\bar{a},\bar{b})$ be the
minimal number $n$ such that for some
$\bar{a}_{0}(=\bar{a}),\bar{a}_{1},...,\bar{a}_{n}(=\bar{b})$
any pair $\bar{a}_{i},\bar{a}_{i+1}$ extends to an infinite
indiscernible sequence.
Newelski has proved in \cite{new}
that a type-definable Lascar strong type has finite diameter
and if the theory is G-compact then there is a finite bound on
the diameters of Lascar strong types.
It is worth noting that in the proof of Theorem \ref{thr}
we explicitely biuld a sequence of Lascar strong type
(of $C_{4i}$'s) with growing finite diameters.
\bigskip

The research was supported by KBN grants 2 P03A 007 19
and 1 P03A 025 28.
The example of a non-G-compact $\aleph_{0}$-categprical theory
was found in 2002 when the author held a visiting position at
Institute of Mathematics of Polish Academy of Sciences.


\end{document}